\documentclass[3p,times]{elsarticle_preprint}

\usepackage{amssymb}

\usepackage{amsmath} 
\usepackage{amsfonts}
\usepackage{textcomp}
\usepackage{mathrsfs}
\usepackage{float}
\usepackage{caption}
\usepackage{subfigure}
\usepackage{mathabx}
\usepackage{graphicx}
\usepackage{verbdef}
\usepackage{mathabx}

\usepackage{listings}
\usepackage{fancyvrb}
\usepackage{caption}
\usepackage{enumitem}
\usepackage{cprotect}
\usepackage{verbatim}
\usepackage{amsmath}
\usepackage{tikz}
\usepackage{color}

\newtheorem{theorem}{Theorem}[section]

\newtheorem{property}[theorem]{Property}

\newtheorem{definition}[theorem]{Definition}

\journal{Applied Mathematics and Computation}

\begin{document}
\newcommand{\NN}{\mathbb{N}}
\newcommand{\ZZ}{\mathbb{Z}}
\newcommand{\QQ}{\mathbb{Q}}
\newcommand{\RR}{\mathbb{R}}
\newcommand{\CC}{\mathbb{C}}
\newcommand{\II}{\mathbb{I}}

\begin{frontmatter}



\title{Hermite-Birkhoff interpolation on scattered data on the sphere and other manifolds}


\author{Giampietro Allasia}
\ead{giampietro.allasia@unito.it}

\author{Roberto Cavoretto}
\ead{roberto.cavoretto@unito.it}

\author{Alessandra De Rossi\corref{cor1}}
\ead{alessandra.derossi@unito.it}
\cortext[cor1]{Corresponding author.}

\address{Department of Mathematics \lq\lq G. Peano\rq\rq, University of Torino, via Carlo Alberto 10, I--10123 Torino, Italy}

\begin{abstract}
The Hermite-Birkhoff interpolation problem  of a  function given on arbitrarily distributed points on the sphere and other manifolds is considered. Each proposed interpolant  is  expressed as a linear combination of  basis functions, the combination coefficients being incomplete Taylor expansions of the interpolated function at the interpolation points. The basis functions  depend on the geodesic distance, are orthonormal with respect to the point-evaluation functionals, and have all  derivatives equal zero up to a certain order at the interpolation points. A remarkable feature of such interpolants, which belong to the class of partition of unity methods, is that their construction  does  not require solving  linear systems. Numerical tests are given to show the interpolation performance.
\end{abstract}



\end{frontmatter}


\section{Introduction}

In  previous papers  we  concerned with Hermite-Birkhoff interpolation of a function given on arbitrarily distributed points on Euclidean spaces \cite{bracco11b} (see also e.g. \cite{dellaccio16a,dellaccio16b}), and with Lagrange and Hermite-Birkhoff interpolation of a function given on arbitrarily distributed points on Banach spaces \cite{bracco11a, bracco11c}. In particular, in \cite{bracco11b} we considered data on a general  domain $\Omega \subset \RR^{d+1},~d\ge 1$, without using any topological information about $\Omega$. Nevertheless  many applications provide us with additional information on the underlying domain. For example, problems coming from geology, geophysics,  metereology, oceanography, satellite-based techniques, etc.,   often relate to the entire earth or a large part of it, so that the unit sphere would be an appropriate model and the additional information would lead to a better approximant. The considered  data typically represent some physical phenomena, such as temperature, rainfall, pressure, and gravitational forces, measured at various points on the surface of the earth, possibly at different times.    

The sphere is a particularly interesting example of a connected compact smooth manifold, even because  considerations developed about interpolation on the sphere can be extended to other manifolds. We think convenient to discuss, as long as possible,  a general framework, though in practice  the most interesting manifolds are smooth two-dimensional manifolds, i.e.  surfaces in $\RR^3$. In fact, a problem that occurs frequently in science and engineering is to recover from a surface in three dimensions a real valued function  that interpolates to a given set of data (see e.g. \cite{Li15, Yao15}).  

The  generalized Hermite  interpolations (or Hermite-Birkhoff interpolation, see e.g. \cite{fass15}, Chapter 19) on scattered data by means of basis functions depending on the distance have been developed since 1992, when appeared the pioneering  paper by Wu  \cite{wu92}. Since then, the interest in this topic seems to have increased significantly (see  e.g. the pertinent chapters  in \cite{fass07, fass15, wendl05}). A number of authors have also considered the  Hermite  interpolation setting on scattered data  on the sphere (see e.g. \cite{freed82, ron96, fass99, narc02, macedo11}) or even general Riemannian manifolds \cite{dyn99, narc95}. 

The point of view we follow in this work is deeply different from that in the just quoted papers.  The matter  is that we do not consider a radial basis function method but a cardinal (radial) basis function method. In this way, to get the Hermite-Birkhoff interpolant we have not to solve a system of linear equations. In fact, the interpolant is directly expressed  as a linear combination of  basis functions, which depend on the geodesic distance, are orthonormal with respect to the point-evaluation functionals, and have all  derivatives equal zero up to a certain order at the interpolation points. The  coefficients of the linear combination are incomplete Taylor expansions of the interpolated function at the interpolation points. On the other hand, our interpolation method is strictly linked up with papers which discuss Lagrange interpolation  by  partition of unity methods, namely    Shepard-like  methods, on the sphere  (see, in particular,  \cite{cavoretto10a, cavoretto10b, cavoretto12, cavoretto14, derossi05, derossi07b}). Our method also shows an interesting analogy with articles by Franssens \cite{franssens04, franssens99}.

This paper is organized as follows. In Section 2 explicit expressions of Hermite-Birkhoff interpolants  on manifolds are given. Since these definitions are based on a suitable class of cardinal basis functions, Section 3 describes a general way to construct such basis functions, which depend on geodesic distances on Riemannian manifolds. Moreover, upper bounds for errors in terms of the fill distance are showed.  In Section 4 some basic considerations are pointed out on geometric objects and analytic tools to deal with.  Beside sketching the general background, the main goal of these considerations is to motivate the strategies to be adopted in numerical computation. Sections 5 and 6 discuss the application of the proposed interpolants to numerical recovering  of functions only known on scattered  data on the sphere and other Riemannian manifolds. The reported numerical tests are restricted to surfaces, which are the most interesting and handy cases, but the adopted point of view is more general. Considering the sphere the well-known expression of the geodesic distance is available, but for other manifolds it is necessary to consider approximations of their  geodesic distances.


\section{Hermite-Birkhoff Interpolation on Manifolds} 

  Now, we  define Hermite-Birkhoff interpolation on Riemannian manifolds, in particular on the higher-dimensional sphere. The manifolds we consider are supposed to enjoy suitable properties, as specified in the following (see Sections 4--6 for geometric details).
\begin{definition}\label{Ba}
\noindent 
Let us consider  a  $m$-dimensional Riemannian manifold ${\cal M}\subset\RR^{d+1}, ~d\ge 2$,   
an open set $U=\{u\equiv (u_1, \ldots, u_{d+1})\in\RR^{d+1}\}\subset {\cal M}$, a function $\varphi:U\to\RR^d$ which maps $U$ homeomorphically to the open set $V=\{v\equiv (v_1, \ldots, v_m)\in\RR^m\}:=\varphi (U)$ so that $v=\varphi (u):=v(u)$ and $u=\varphi^{-1}(v):=u(v)$.
Let ${\cal X}=\lbrace  z_1,\ldots, z_n\rbrace\subset U$ be a set of  distinct points, possibly scattered,  with associated finite sets $\Delta_1,\ldots,\Delta_n\subset \NN_0^m$. The {\rm Hermite-Birkhoff interpolation problem} from  $U$ to $\RR$  consists in finding a function  $H:U\rightarrow \RR$ which satisfies the interpolation conditions
\begin{equation}
D^\beta H( z_i):=\frac{\partial^{|\beta|}H( z_i)}{\partial v_1^{\beta_1}\cdots \partial v_m^{\beta_m}}=f_{i \beta}, \qquad \beta\in\Delta_i,~ 
\quad i=1,\ldots,n,
\label{ba}
\end{equation}
where $\beta=(\beta_1, \ldots, \beta_m)$, $\vert\beta\vert=\beta_1+\cdots +\beta_m$, and the $f_{i \beta}\in \RR$ are  given values to be interpolated.
It is  assumed that $H\in C^k(U)$, where $k=\max\lbrace\vert \beta\vert: \beta\in\Delta_i,\,\hbox{for some }i,\, 1\le i\le n\rbrace$.
\end{definition}
\noindent In the following, we will also think that the $f_{i\beta}$ are values assumed by an underlying  function $f:U\rightarrow\RR$, $f\in{\cal C}^k(U)$, so that the conditions \eqref{ba} take the form
\begin{equation}
D^\beta H( z_i)=D^\beta f( z_i), \qquad \beta\in\Delta_i, \quad i=1,\ldots,n.
\label{bb}
\end{equation}
In general, the values of $f$ and of some its derivatives are known only at the points of ${\cal X}$.

A constructive solution to the interpolation problem \eqref{ba} can be given by introducing a suitable class of {\sl cardinal basis functions}, which can be   defined as follows.
\begin{definition}\label{Bb}
\noindent Given a set of distinct points ${\cal X}=\{ z_i,~1\le i\le n \}$, arbitrarily distributed in the open set  $U\subset {\cal M}$, the functions $g_j: U\to\RR,~1\le j\le n$, are cardinal basis functions with respect to ${\cal X}$ if they satisfy for all $u\in U$ the conditions
$$
g_j\in C^k(U),\qquad g_j(u)\ge 0  ,\qquad \sum_{j=1}^n g_j(u)=1 ,\qquad g_j( z_i)=\delta_{ji} ,
$$
where  
$\delta_{ji}$ is the Kronecker delta, and also satisfy the additional property
\begin{equation}
D^\beta g_j(z_i)=0, \qquad \beta\in\Delta_i,~\vert \beta\vert\ne 0, \quad i=1,\ldots,n. 
\label{bc}
\end{equation}
\end{definition}
\noindent It is clear that an interpolant based on these weights must be considered as a {\sl partition of unity method}.

It can be easily checked that the following result holds.
\begin{property}\label{Bc}
The interpolation conditions \eqref{bb} are satisfied  by the interpolant
\begin{align}
H(u)=\sum_{i=1}^n ~T\bigl(u;f, z_i,\Delta_i)~g_i(u),\label{bd}
\end{align}
where
\begin{align*}
T\bigl(u;f,  z_i,\Delta_i):=\sum_{\beta\in\Delta_i} \frac{D^\beta f(z_i)}{\beta_1!\cdots\beta_m!}  \big(v-v(z_i)\big)^\beta=\sum_{\beta\in\Delta_i} \frac{D^\beta f( z_i)}{\beta_1!\cdots\beta_m!} \big(v_1-v_1(z_i)\big)^{\beta_1}\cdots \big(v_m-v_m(z_i)\big)^{\beta_m},
\nonumber
\end{align*}
is formally an incomplete Taylor expansion of $f$ at $ z_i$, in the sense that it only includes the partial derivatives whose orders belong to $\Delta_i$. The interpolant \eqref{bd} can also be expressed in the form
\begin{equation}
H(u)=\sum_{i=1}^n \sum_{\beta\in\Delta_i} ~D^\beta f(z_i)~g_{i\beta}(u),
\label{be}
\end{equation}
where
\begin{equation}
g_{i\beta}(u):=
\frac{\big(v-v(z_i)\big)^\beta}{\beta_1!\cdots\beta_m!}g_i(u)=
\frac{\big(v_1-v_1(z_i)\big)^{\beta_1}\cdots \big(v_m-v_m(z_i)\big)^{\beta_m}}{\beta_1!\cdots\beta_m!}g_i(u),  
\label{bf}
\end{equation}
with $\beta\in\Delta_i,~ 1\le i\le n$.
\end{property}
\noindent The formula \eqref{be} highlights that the interpolant is essentially constructed by using the $g_{i\beta}$  as basis functions.

\medskip
 
The interpolant \eqref{bd} enjoys the usual properties of cardinal basis interpolants.
\begin{property}\label{Bd}
\noindent There hold  the following inequalities:
\begin{align*}
& a)\quad \Vert H(u)\Vert\le \max_i \Vert T\bigl(u;f, z_i,\Delta_i)\Vert,\\
& b)\quad \Vert f(u)-H(u)\Vert\le \sum_{i=1}^n g_i(u)\Vert f(u)-T\bigl(u;f,  z_i,\Delta_i)\Vert\le \max_i \Vert f( u)-T\bigl( u;f,  z_i,\Delta_i)\Vert,
\end{align*}
where the $i$-th term in the sum may be interpreted as the local error at the point $ z_i$.
\end{property}


\section{Cardinal Basis Functions on Manifolds}

A classical way to construct cardinal basis functions defined on $\RR^{d+1}$ is Cheney's method (see \cite{cheney} and \cite{cheney00}, pp. 67--68), which can be used for manifolds as well, if we adopt a suitable  distance.  

\begin{theorem}\label{Ca}
Let us consider $U\subset {\cal M}$ as in Definition \ref {Ba} and let $\alpha:U\times U \rightarrow \RR^+_0$ be a continuous and bounded function, such that $\alpha(u,z_i)>0$ for all $u\in U$, $u\ne z_i$, and $\alpha (z_i, z_i)=0$ for all $z_i\in {\cal X}$. Moreover, let each $\alpha(u,z_i)$ be $k$-times continuously differentiable on $U$ such that
\begin{equation*}
[D^{\beta} \alpha(u,z_i)]_{u=z_i}=0, \qquad i=1,\ldots,n, \qquad 0<\vert\beta\vert\le k. \notag
\end{equation*}
The corresponding {\rm cardinal basis functions}
\begin{equation}
g_i(u)={{\displaystyle{\prod_{j=1, j\neq i}^n \alpha (u,z_j)}}
\over{{\displaystyle\sum_{k=1}^n\prod_{j=1, j\neq k}^n \alpha (u,z_j)}}},\qquad i=1,\ldots,n,  
\label{ca}
\end{equation}
are continuous and satisfy
\begin{equation}
D^\beta g_i(z_j)=0, \qquad 0<\vert\beta\vert\le k, \quad i,j=1,\ldots,n. 
\label{cb}
\end{equation}
\end{theorem}

\smallskip
\noindent
{\bf Proof:} This result is essentially the $d-$dimensional case of the main theorem in \cite{bracco11b}. 
\hfill $\square$

\medskip
\noindent
 The $g_i$ in \eqref{ca}  can also be represented in the {\sl barycentric form}
$$
g_i(u)=\frac{\displaystyle{1/\alpha(u,z_i)}}{\displaystyle{\sum_{k=1}^n 1/\alpha(u,z_k)}},\qquad g_i(z_i)=1,\qquad i=1,\ldots,n,
$$
which is more convenient from a computational point of view \cite{allasia97}.

A natural choice is defining $\alpha$ using the distance between points. Since we are considering points on the manifold ${\cal M}$, we   take  the geodesic distance $d_g$ 
(see Definition \ref{De} below)  
and  define  $\alpha$ in the general form
\begin{equation}
\alpha(u, w)=\vartheta (d_g(u, w))  ,
\label{cd} 
\end{equation}
which obviously must satisfy the assumptions of Theorem \ref{Ca}.  
In particular, in (\ref{cd}) we  may consider the choice
\begin{equation}
\alpha(u,w)= (d_g(u,w))^\mu, \qquad  \mu\in\RR^+ ,~ \mu\ge k,  ~u,w\in U,
\label{ce} 
\end{equation}
 which ensures both the vanishing of the derivatives at the nodes and the regularity assumptions, and among the  possible choices is the most direct.  
Other interesting choices are  (see e.g. \cite{cavoretto10a}):
\begin{align*} 
&\alpha_\gamma(u,w)=\frac{\exp(\gamma\hskip 1pt (d_g(u,w))^\mu)}{(d_g(u,w))^\mu}, \quad \gamma>0,~ \mu\ge k,
\nonumber
\cr
&\alpha_\delta (x,y)={\exp(\delta\hskip 1pt d_g(u,w)^\mu)}, \quad \delta \ge 0,~\mu\ge k,
\end{align*} 
both of them being rapidly decreasing.

As a result  of the choice (\ref{ce}), we obtain the cardinal basis functions
\begin{equation}
g_i(u)=\frac{(d_g(u, z_i))^{-\mu}}{\sum_{k=1}^n (d(u, z_k))^{-\mu}}, \quad i=i, \ldots, n,
\label{cf}
\end{equation}
but, for  computational reasons, in many cases it may be  preferable to use a localized version of the cardinal basis functions \eqref{cf}, that is,
\begin{equation}
\tilde g_i(u)=\frac{\tau_i(u)(d_g(u, z_i))^{-\mu}}{\sum_{k=1}^n \tau_k(u)(d_g(u, z_k))^{-\mu}}, \label{cg}
\end{equation}
where $\tau_i:U\to\RR^+_0, ~\tau_i\in{ C}^k(U)$, such that
\begin{equation}
\tau_i(u)=
\begin{cases}
>0, ~\hbox{for } u:d_g(u,z_i)<\delta,\\ 
=0, ~\hbox{for } u:d_g(u,z_i)\ge\delta, \\
\end{cases}
\label{cga} 
\end{equation}
and  $\delta>0$ is a suitably chosen value. 
Hence, each function $\tilde g_i$  vanishes outside the neighborhood of $z_i$  such that   $d_g(u,z_i)\ge\delta,~u\in U$. It can be easily proved that the functions $\tilde g_i$, $i=1,\ldots,n$, are cardinal and enjoy the vanishing property on derivatives \eqref{bc}. 
The localization can be obtained by taking, for instance,
$$
\tau_i(u)=\Big(1-\frac{d_g(u, z_i)}{\delta}\Big)^{k+1}_+,\quad u\in U, \quad i=1,\ldots,n.
$$
A simpler, but a little naive, choice is 
\begin{equation}
\tau_i(u)=
\begin{cases}
1, ~\hbox{for } u:d_g(u,z_i)<\delta,\\ 
0, ~\hbox{for } u:d_g(u,z_i)\ge\delta. \\
\end{cases}
\nonumber 
\end{equation}

For the  Hermite-Birkhoff interpolant with cardinal basis functions \eqref{cg}
\begin{align}
\tilde H(u)=\sum_{i=1}^n ~T\bigl(u;f,z_i,\Delta_i)~ \tilde g_i(u),
\label{ch}
\end{align}
we can give more significant error estimates than for the basic case \eqref{bd}.
Let
$q\in\NN $ be defined such that each Taylor-type expansion $T(u; f, z_j, \Delta_j)$ is a complete Taylor expansion up to order $q$, plus other terms of higher degree.
For any $f: U\to\RR$ with $f\in{ C}^q(U)$ and for any $u\in U$, we have, since the cardinal basis functions $\tilde g_i$ are a partition of unity,
\begin{align}
&\vert f(u)-\tilde H(u)\vert=\Bigg\vert\sum_{i=1}^n~f(u)\,\tilde g_{i}(u)-\sum_{i=1}^n~T(u; f, z_i, \Delta_i)\,\,\tilde g_{i}(u)\Bigg\vert=\notag \\
&\Bigg\vert \sum_{i=1}^n~\big[f(u)-T(u; f, z_i, \Delta_i)\big]\,\,\tilde g_{i}(u)\Bigg\vert\le \sum_{i=1}^n~\big\vert f(u)-T(u; f, z_i, \Delta_i)\big\vert\,\,\tilde g_{i}(u), 
\label{ck}
\end{align}
each $\tilde g_i$ being  non-zero only inside the ball of radius $\delta$ centered at $z_i$. Now, since each $T(u;f,z_i,\Delta_i)$ is a Taylor expansion complete up to order  $q$, we can use the estimate 
\begin{equation}
\vert f(u)-T(u; f, z_i, \Delta_i)\vert \le C_i \Vert v-v(z_i)\Vert^{q+1}, 
\label{cl}
\end{equation}
where $C_i\in\RR^+$ is a suitable constant and  $\Vert\cdot\Vert$ is the Euclidean norm.
Since $\Vert v-v(z_i)\Vert$ is less than or equal to the geodesic distance $d_g(u, z_i)$, 
it  follows
\begin{equation}
\vert f(u)-T(u; f, z_i, \Delta_i)\vert \le C_i d_g^{q+1}(u,z_i). 
\label{cm}
\end{equation}
Inserting \eqref{cm} in \eqref{ck} and exploiting again the partition of unity property, since $\vert d_g(u,z_i)\vert<\delta$, we get
\begin{equation*}
\vert f(u)-\tilde H(u)\vert\le \sum_{i=1}^n C_i d_g^{q+1}(u,z_i)\tilde g_i(u)
\le \delta^{q+1}\sum_{i=1}^n C_i \tilde g_i(u)\le C\delta^{q+1}, 
\end{equation*}
with $ C=\max_i C_i$.
Moreover, if we set the localization radius $\delta=K h_{U,{\cal X}}$, where $h_{U,{\cal X}}$ is the so-called {\sl fill distance}, that is,
\begin{align}
h_{U,{\cal X}}:=\sup_{u\in U }\inf_{z_i\in {\cal X}} d_g(u,z_i).
\label{cra}
\end{align}
and $K\ge 1$, we obtain the estimate
\begin{equation}
\vert f(u)-\tilde H(u)\vert\le CK h_{U,{\cal X}}^{q+1}. 
\label{cn}
\end{equation}


\section{Some Facts from Geometry}

Considering the Hermite-Birkhoff interpolation  problem on the sphere $\mathbb{S}^{d}$  and   Riemannian manifolds in $\RR^{d+1}$, some basic considerations are to be pointed out on the  geometric objects and analytic tools to deal with. In fact, besides giving the general background, they  affect  deeply  numerical computation strategies.

 Let us consider first the interpolation problem on the sphere. It is well known that the spherical earth model is navigated using flat maps, collected in an atlas, and no single flat map can represent the entire earth. Similarly, the sphere $\mathbb{S}^d\subset\RR^{d+1}$ can be described using an atlas of  charts, each chart mapping part of the sphere to a subset of $\RR^d$.   Precisely, for every $u\in \mathbb{S}^d$ there exist an open set $U\subset \mathbb{S}^d$ with $u\in U$
and a mapping $\varphi:U\to\RR^{d}$ that maps $U$ homeomorphically to the open set $V:=\varphi(U)$.  The pair $(U, \varphi)$ is  a  {\rm chart}  of $U$ and a  collection ${\cal A}=\{(U_\alpha, \varphi_\alpha)\}$ of charts, which covers the sphere,  is an {\rm atlas}. 
Charts in an atlas may overlap and a single point of the sphere may be represented in different charts. 
Given two  overlaping charts $(U_\alpha, \varphi_\alpha)$ and $(U_\beta, \varphi_\beta)$, a transition function, that is, a coordinate transformation,   can be defined which goes from $\varphi_\alpha(U_\alpha\cap U_\beta)\subset\RR^d$ to the sphere and then back to $\varphi_\beta(U_\alpha\cap U_\beta)\subset\RR^d$. 
A chart $(U, \varphi)$ is of class $C^k$ if $\varphi^{-1}\in C^k(\varphi(U))$, whereas a  $C^k$-atlas consists of $C^k$-charts and $C^k$-transition functions.

 Through the chart $(U, \varphi)$  the neighborhood $U$ inherits the  coordinates given on the Euclidean space $\RR^d$ and  the homeomorphism $\varphi$ leads us to describe $U$ as a locally Euclidean space. In fact, considering $\varphi^{-1}:V\to U$ we have that 
the coordinates $u_1, \ldots, u_{d+1}$ of a point  $u\in U$ can be given by  $d+1$ parametric equations  
\begin{equation}
u_1=u_1(v_1, \ldots, v_{d}),
\ldots ,
u_{d+1}=u_{d+1}(v_1, \ldots, v_{d}),
\label{da} 
\end{equation}
where  the parameters  $v_1, \ldots, v_d$  identify the  point    $v=\varphi(u)\in V\subset \RR^d$. 
Then, the map $\varphi^{-1}$ can be written in terms of its components as $\varphi^{-1}(v_1, \ldots, v_d)=\big(u_1(v_1, \ldots, v_d), \ldots, $ $u_{d+1}(v_1, \ldots, v_d)\big)$.

An atlas is not unique as the sphere can be covered in multiple ways using different combinations of charts. To  describe a possible atlas for the sphere ${ S}^{d}$,   we consider for any $u^*\in { S}^{d}$ the open neighborhood of $u^*$ in ${ S}^{d}$ given by $U^+_{u^*}=\{u\in { S}^{d} : (u,u^*)>0\}$, where $(\cdot, \cdot)$ is the usual inner product in $\RR^{d+1}$. Choosing the coordinate system in $\RR^{d+1}$ so that the vector $u^*$ has, to say,  components $(0,0,\ldots,1)$ we have that the neighborhood $U^+_{u^*}$ can be  homeomorphically projected on an open set  in $\RR^d$. Similarly, we may consider the neighborhood $U^-_{u^*}=\{u\in { S}^{d} : (u,u^*)<0\}$. Since ${ S}^{d}$ can be thought  as the union of a suitable number of charts, we get an atlas   and the considered parametrization holds. In particular, considering the two-dimensional sphere $\mathbb{S}^2$   an atlas of six charts is obtained which covers the entire sphere. Otherwise for the sphere $\mathbb{S}^2$, as well as for $\mathbb{S}^d$, it may be sometimes convenient to use  spherical  coordinates (see e.g. \cite{fass99, narc95}).   Choosing one or another atlas has significant effects, especially for  actual applications and their numerical treatment. 

 Referring to a real function $\psi$ defined on the sphere $\mathbb{S}^d$,  
we say that $\psi$ is $k$-times differentiable on  $\mathbb{S}^d$, or $\psi\in C^k(\mathbb{S}^d)$, if $\psi\circ\varphi^{-1}\in C^k(\varphi(U))$ for every chart $(U, \varphi)$ of ${\mathbb{S}^d}$.  
The function $\psi:{\mathbb{S}^d}\to\RR$   inherits the local coordinates $v_1, \ldots, v_{d}$ of the chart $U$, since $\psi(u)=\psi\circ \varphi^{-1}(v)$.

The considerations on the sphere, just seen,   can be extended to other manifolds. It is useful  to recall a formal definition of  a topological manifold (see e.g. \cite{wendl05, boothby75}).
\begin{definition}\label{Da}
\noindent 
A set ${\cal M}\subset \RR^{d+1}$ is called a {\rm topological manifold} of dimension $m$,
if it is a Hausdorff space with a countable  basis of open sets such that for every $u\in {\cal M}$ there exist an open set $U\subset {\cal M}$ with $u\in U$   and a mapping $\varphi:U\to\RR^{m}$  that  maps $U$ homeomorphically to the open set $V:=\varphi(U)\subset \RR^m$. The pair $(U, \varphi)$ is called 
a coordinate neighborhood of $u$ or 
a {\rm chart} and for every $u\in U$ the vector $\varphi (u)=(v_1(u), \ldots, v_m(u))\in \RR^m$ represents the {\rm local coordinates} of $u$ in $V$. A chart is of class $C^k$ if $	\varphi^{-1}\in C^k(\varphi(U))$. A collection ${\cal A}=\{(U_\alpha, \varphi_\alpha)\}$ of $C^k-$charts is called a $C^k-${\rm atlas} of ${\cal M}$ if
 the sets $U_\alpha$ cover ${\cal M}$ and, moreover,
for any $U_\alpha, U_\beta$ with $U_\alpha\cap U_\beta\ne\emptyset$ the {\rm transition maps} $\varphi_\beta \circ\varphi_\alpha^{-1}$ and $\varphi_\alpha \circ\varphi_\beta^{-1}$ are in $C^k$ on $\varphi_\alpha(U_\alpha\cap U_\beta)$ and $\varphi_\beta(U_\alpha\cap U_\beta)$, respectively.
Finally, a manifold ${\cal M}$ is called a $C^k-${\rm manifold} if it possesses a $C^k-$atlas.
\end{definition}

The smoothness of a function $f:{\cal M}\to \RR$ is defined by the smoothness of $f\circ\varphi^{-1}$ with a chart $(U, \varphi)$, as it is pointed out in the following: 
\begin{definition}\label{Df}
We say  that  $f: {\cal M}\to\RR$  is $k$-times differentiable on ${\cal M}$ 
or $f\in C^k({\cal M})$ provided that for every chart $(U, \varphi)$ of  ${\cal M}$
the composition
$$
f\circ \varphi ^{-1}: \varphi (U)\to \RR
$$ 
is $k$-times differentiable.
\end{definition}

\noindent It is important to realize that the definition of differentiability of a real-valued function on a $C^k$-manifold does not depend on the choice of the chart.

To introduce on a $C^k$-manifold  ${\cal M}\subset\RR^{d+1}$, $k\ge 1$, the notions of  length and distance,  each tangent space must be equipped with an inner product, so that it varies smoothly from point to point. 
The tangent space $T_u{\cal M}$ for a point $u\in{\cal M}$  is the space formed by the tangent vectors to all the curves in ${\cal M}$ passing  through $u$. 
Here, a vector $\tau$ is a tangent vector in $u\in{\cal M}$ if there exists a differentiable curve $\gamma(t)$ on ${\cal M}$, depending on a parameter $t$ with $-\varepsilon\le t\le  \varepsilon$  for some $\varepsilon>0$, such that $\gamma(0)=u$ and $\gamma'(0)=\tau$, where
$$
\gamma'(0):=\frac{d}{dt}\varphi\circ\gamma(t)\Big|_{t=0}.
$$
It turns out that $T_u({\cal M})$ is a $m-$dimensional vector suspace of $\RR^{d+1}$ and that a basis is given by
\begin{equation}
\frac{\partial \varphi^{-1}}{\partial v_1}(\varphi(u)), \ldots, \frac{\partial \varphi^{-1}}{\partial v_m}(\varphi(u)).
\label{db}
\end{equation}
It is interesting to note that the tangent space can be thought as  the best linear approximation to ${\cal M}$ in $u$.

 More explicity, let us consider a chart $(U, \varphi)$ of a manifold ${\cal M}$ with  parametric equations 
\begin{equation}
u_1=u_1(v_1, \ldots, v_m), \ldots, u_{d+1}=u_{d+1}(v_1, \ldots, v_m),
\nonumber
\end{equation}
and a curve $\gamma:[-\epsilon, +\epsilon]\to{\cal M}$ 
on $U$ whose  equations in $\varphi(U)$ are
\begin{equation}
v_1=v_1(t), \ldots, v_{m}=v_{m}(t).
\nonumber
\end{equation}
Substituting the latter equations into those of the chart, we get the  equations of the curve $\gamma(t)$ on $U$ as a function of $t$, that is, $\gamma(t)=u(v(t))$. Then, differentiating we obtain  the
tangent vector to the curve  at the point $u=\gamma(0)$ 
\begin{equation}
\gamma'(0)
=\frac{\partial u}{\partial v_1} \frac{dv_1}{dt}+ \ldots+\frac{\partial u}{\partial v_m}\frac{dv_m}{dt}\Big|_{t=0}
\label{de}
\end{equation}
and the  basis vectors  are 
$$
\frac{\partial u}{\partial v_1}=\Big(\frac{\partial u_1}{\partial v_1},\ldots, \frac{\partial u_{d+1}}{\partial v_1}\Big), \ldots, \frac{\partial u}{\partial v_m}=\Big(\frac{\partial u_1}{\partial v_m},\ldots,\frac{\partial u_{d+1}}{\partial v_m}\Big),
$$
namely  (\ref{db}).

 To operationalize the Hermite-Birkhoff interpolation technique  we  start  considering the geometric problem involving the computation of lengths of curves lying on a surface ${\cal M}\subset\RR^3$. The key idea is essentially based on replacing an infinitesimal element of a smooth curve by the corresponding element of its tangent plane.
As a significant  example, which concerns the sphere and other major surfaces, let us consider a  surface ${\cal M}$, a curve $\gamma$ on it and a point $u=(u_1(v_1,  v_2), ~u_2(v_1, v_2), u_3(v_1,  v_2))$ on $\gamma$.
Taking the arclength $s$ of the curve as  a parameter, 
the vector $du/ds$ 
is of unit length and  we have
from \eqref{de} for $m=2$
$$
ds^2=\Big(\frac {\partial u}{\partial v_1},\frac {\partial u}{\partial v_1}\Big)dv_1^2+2\Big( \frac {\partial u}{\partial v_1},\frac {\partial u}{\partial v_2}\Big) dv_1dv_2+\Big( \frac {\partial u}{\partial v_2},\frac {\partial u}{\partial v_2}\Big)dv_2^2,
$$
where $(\cdot, \cdot)$ is the scalar product.
Making use of the notations
$$
g_{11}=\Big( \frac {\partial u}{\partial v_1},\frac {\partial u}{\partial v_1}\Big), \quad g_{12}=g_{21}=\Big(\frac {\partial u}{\partial v_1},\frac {\partial u}{\partial v_2}\Big), \quad g_{22}=\Big( \frac {\partial u}{\partial v_2},\frac {\partial u}{\partial v_2}\Big),
$$
we obtain the {\sl first fundamental quadratic form} of the surface
$$
ds^2=g_{11}dv_1^2+2g_{12}dv_1dv_2+g_{22}dv_2^2.
$$
The components $g_{ij},(i,j=1,2),$ of the metric  form the entries of a $2\times 2$ symmetric matrix, namely $ds^2$ is a positive quadratic form related to  this matrix.
The first fundamental quadratic form of a surface provides the expression for the length $ds$ of an infinitesimal arc and the length of a finite curve lying on the surface is obtained from it by integration. More precisely, if a curve on the surface is given by the equation $\gamma(t)=u(v(t)),
~t_1\le t\le t_2,$ its length is
$$
L(\gamma; t_1, t_2)=\int_{t_1}^{t_2}\Big[g_{11}\Big(\frac{dv_1}{dt}\Big)^2+2g_{12}\frac{dv_1}{dt}\frac{dv_2}{dt}
+g_{22}\Big(\frac{dv_2}{dt}\Big)^2\Big]^{1/2}dt.
$$

  To handle this idea in a more general situation,  we  recall the concept of Riemannian manifold.
\begin{definition}\label{Dc}
\noindent  A $C^k$-manifold is called a $C^k-${\rm Riemannian manifold} if for every $u\in{\cal M}$ there exists  a positive definite inner product 
$g_u:T_u({\cal M})\times T_u({\cal M})\to \RR$ such that for every chart $(U, \varphi)$ the $m^2$ functions
\begin{equation}
g_{ij}^\varphi(u):=g_u\Big( \frac{\partial \varphi^{-1}}{\partial v_i}(\varphi(u)),  \frac{\partial \varphi^{-1}}{\partial v_j}(\varphi(u))\Big), \quad i,j=1, \ldots, m,
\label{df}
\end{equation}
are in $C^k(V)$ with $  V=\varphi(U)$.
The family of   $g_{ij}^\varphi(u)$, assuming compatibility among different charts, is called  a {\rm Riemannian metric} on 
${\cal M}$, the $g_{ij}^\varphi(u)$ are the component of the metric and form the entries of a $m\times m$ symmetric matrix. The first fundamental quadratic form associated to the metric is
\begin{equation}
ds^2=\sum_{i, j=1}^m g^\varphi_{ij}(u)dv_idv_j.
\nonumber
\end{equation}
\end{definition}

\noindent All differentiable manifolds (of constant dimension) can be given the structure of a Riemannian manifolds. The Euclidean space itself carries a natural structure of Riemannian manifold, where the tangent spaces are naturally identified with the Euclidean space itself and  the scalar product of the space is  the standard scalar product. Precisely, with ${\partial u}/{\partial v_i}$ identified with the $i$-th standard basis vector $e_i = (0, \ldots, 1, \ldots, 0)$, the  (canonical) Euclidean metric over an open subset $U\subset \RR^{d+1}$ is defined by
$ g^{\mathrm{can}}_{ij}=(e_i,e_j) = \delta_{ij}$.
Many familiar curves and surfaces, including for example all $d-$spheres, are specified as subspaces of a Euclidean space and inherit a metric from their imbedding.

Finally, we use the Riemannian metric to define the length of a curve on ${\cal M}$.
\begin{definition}\label{De}
\noindent  Suppose that ${\cal M}$ is a  connected $C^k-$Riemannian manifold. Let $u, w\in {\cal M}$ be two distinct points and let $\gamma:[a, b]\to {\cal M}$ be a piecewise $C^1$ curve that connects these points, i.e., $\gamma(a)=u, \gamma(b)=w$. Then the length of $\gamma$ is  expressible in one of the equivalent forms
\begin{align}
L(\gamma;a,b):=&\int_a^b \Vert \gamma'(t) \Vert dt=\int_a^b \Big[g_{\gamma(t)}\Big(\frac{d\gamma}{dt}(t), \frac{d\gamma}{dt}(t)\Big) \Big]^{1/2}dt
\cr
=&\int_a^b\Big(\sum_{i, j=1}^m g^\varphi_{ij}(u)\frac{dv_i}{dt} \frac{dv_j}{dt}\Big)^{{1/2}}dt,
\label{dg}
\end{align}
where the first integrand represents the length of an infinitesimal arc in the Riemannian metric and $\Vert\cdot\Vert$ denotes the norm induced by the inner product on the tangent space.
Supposing ${\cal M}$ to be compact, we set $d_g(u, w)$ to be the infimum over the length of all such curves connecting $u$ and $w$. The shortest of such curves is called the shortest path for $u$ and $w$, and $d_g(u, w)$ is their {\rm geodesic} or {\rm Riemannian distance}.
\end{definition}

\noindent
  If ${\cal M}=\RR^m$ and if $g_u$ is the canonical inner product on $\RR^m$ then our definition of the length of a curve in the Riemannian metric coincides with the classical definition. In this case, $d_g(u, w)=\Vert u-w\Vert_2$, i.e. the shortest curve between two points in $\RR^m$ is the straight line between them.  On the sphere, our definition of geodesic  coincides with the old one, since both denote the length of the shorter portion of the great circle connecting the two points.


\section{Numerical Results on the Sphere}
 In this section we discuss numerical calculation  of  Hermite-Birkhoff interpolation on the sphere. Referring to the framework in  Sections 2 and 3, we develop our considerations   on  $\mathbb{S}^d, ~d\ge 2$,  as  long as possible.

\subsection{Computation of the Geodesic Distance on the Sphere}
 
The $d-$dimensional sphere  represents a case where the concept of atlas is essential.  As described in Section 4, 
suitable charts are given for $\mathbb{S}^d$ by the function $\varphi^+_i(u_1, \ldots, u_{d+1})=(v_1, \ldots, v_{i-1}, v_{i+1},  \ldots, v_d)$,   which projects the subset of $\mathbb{S}^d$ with $u_i>0$ in the subspace $\RR^d$, and by the similar function $\varphi^-_i$ 
with $u_i<0$, being $i=1, \ldots, d+1$. The family of all these charts forms an atlas.

Referring to Definition \ref {De},     the geodesic distance of   any two points $u$ and $w$ of the unit sphere $\mathbb{S}^{d}$, that is, the sphere of unit radius and  centered at the origin, is
\begin{equation}
d_g(u, w)=\arccos(u, w).
\nonumber
\end{equation} 
The geodesic distance $d_g(u,w)$ is always expressible in terms of the Euclidean distance $d_E(u,w)$
in $\RR^{d+1}$.
In fact, we have 
\begin{equation}
d_g(u, w)=
2\arcsin\frac{d_E(u, w)}{2},
\nonumber
\end{equation} 
and, conversely, 
\begin{equation}
d_E(u, w)=\sqrt{2-2(u,w)}=2\sin\frac{d_g(u,w)}{2}.
\nonumber
\end{equation}
 It follows from  the asymptotic expansion of $\arcsin x, ~x\in\RR$, 
\begin{align*}
d_g(u, w)=d_E(u, w)+\frac{1}{24}d\hskip 1pt_E^3(u, w)+O(d\hskip 1pt_E^5(u, w)),
\end{align*}
that is, the difference between  $d_g(u, w)$ and  $d_E(u, w)$ may be very  small if $u$ and $w$ are sufficiently close. {Therefore, using our local interpolation method, the Euclidean distance may be considered as a good approximation of the geodesic distance.}

In general, taking  a  radial function $\varphi(d_E(u,w))$ we obtain a {\it zonal basis function}
$\psi(d_g(u, w))$ with $u,w\in \mathbb{S}^{d}$ by setting $$
\varphi (d_E(u,w))=\varphi(\sqrt{2-2(u,w)})=\varphi(\sqrt{2-2\cos (d_g(u,w)})=\psi(d_g(u, w)).
$$
In particular,  referring to the construction of  cardinal basis functions in Section 3, we remark that the expressions of $\alpha$ in terms of Euclidean and geodesic distances are closely related, because  from (\ref{cd}) we have for suitable functions $\vartheta$ and $\eta$
\begin{equation}
\alpha(u, w)=\vartheta(d_g(u,w))= \vartheta\bigg(2\arcsin\frac{d_E(u,w)}{2}\bigg)=\eta(d_E(u,w)), 
\qquad  u,w\in U.
\label{cda} 
\end{equation}
Hence, known expressions of $\alpha$ in terms of the Euclidean distance can be used as well to get expressions of $\alpha$ in terms of the geodesic distance (see e.g. \cite{hubbert01}).

\subsection{Computation of the Interpolant on the Sphere}

Implementing the interpolation formula $H:U\subset \mathbb{S}^d\to\RR$ in \eqref{bd}, a  problem is to optimize the nearest neighbour searching procedure for spherical data, that is finding in a convenient way the set of  points of ${\cal X}$ closest to any point $u\in U$. It is possible to  use a cell-technique, which consists in a space decomposition into hypercubic cells by overlaying a spatial grid on the sphere. The procedure, which  is based on the optimized Renka's algorithm  for trivariate interpolation \cite{renka88, cavoretto15, cavoretto16},  has been successfully tested for $\mathbb{S}^2$ \cite{derossi05}.  Referring only to $\mathbb{S}^2$, another procedure, as well successfully tested, makes a decomposition of $U\subset{\cal M}$ in strips or spherical zones \cite {allasia11,cavoretto10b}.

Generally speaking, it is convenient to use the  localized version  (\ref{cg}) of the cardinal basis functions, so that more distant points have less influence. Alternatively,  one could  consider exponential-type weights, as in (\ref {cf}), which are strongly decaying as distance increases (see  e.g. \cite{buhmann03}, p. 46). The drawback is that at least one parameter is necessary for the localization and this implies the requirement of determining its optimal value.  The choice of an appropriate value for the localization parameter $\delta$  in  \eqref{cga} determines the efficiency of the local scheme and is  a nontrivial problem. In practice, the localization can be obtained using for  interpolation only the nodes that belong to a convenient neighborhood of the point $u$ considered, i.e. the nodes $z_k$ whose distance from $u$ is smaller than $\delta$.

\begin{figure}[ht!]
		\centering
		\includegraphics[height=.31\textheight]{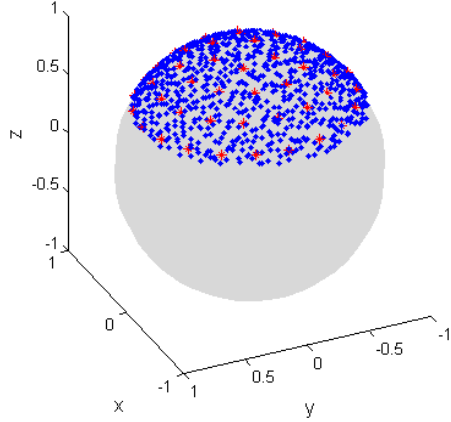}
		\caption{Example of interpolation (blue dot) and evaluation (red star) points on the sphere.}
		\label{fig_sphere}
	\end{figure}

To test the performance of our interpolation method, we must consider in general a {\sl uniform distribution of nodes on the sphere}. To generate uniformly (or quasi-uniformly) distributed random (or pseudo-random) points on  the high-dimensional unit sphere  one can use, in principle, anyone of the suitable algorithms proposed in the literature, but they  are not quite equivalent and only some of them work for $d\ge 3$. The most efficient and fast algorithm is based on the fact that the normal distribution function for a point $u\in\RR^{d+1}$ has a density that depends only on the distance of the point from the origin, so that the points of $\mathbb{S}^d$ have the uniform distribution (see e.g. \cite{knuth96,marsaglia72,muller59,weisstein1,weisstein2}). Since our interpolation algorithm works locally and on suitable charts, we experimented another way to find interpolation points on the unit sphere $\mathbb{S}^2$ $\subset$ $\mathbb{R}^3$ (and in the following on manifolds), considering $n$ Halton points \cite{wong97} on the spherical cap of $\mathbb{S}^2$ for $z>0.5$. Interpolation errors are instead evaluated on a nearly uniform distribution of $n_{eval}=50$ spiral points belonging to the considered chart (see e.g. \cite{saff97,cavoretto10b}). An example of interpolation and evaluation points defined on the chart of $\mathbb{S}^2$ is shown in Figure \ref{fig_sphere}. 

Therefore, to investigate accuracy of the Hermite-Birkhoff interpolant, we compute the Maximum Absolute Error (MAE) and the Root Mean Square Error (RMSE) given by
\begin{align} \label{MAE_RMSE}
	{\rm MAE} = \max_{1\leq i\leq n_{eval}} |f(u_i) - {\tilde{H}}(u_i)|, \qquad {\rm RMSE} = \sqrt{\frac{1}{n_{eval}}\sum_{i=1}^{n_{eval}} |f(u_i) - {\tilde{H}}(u_i)|^2}.
\end{align}

In order to get an idea of the  distribution of points in the set ${\cal X}\subset U$ (see Definition \ref {Ba}) and, in particular, of their uniformity and density, we consider two common indicators, that is,   
the {\sl separation distance}
\begin{equation}
q_{U,{\cal X}}:=\frac{1}{2}\min_{i\ne j}d_g(z_i, z_j),
\label{cdb} 
\end{equation}
and the {\sl fill distance} \eqref{cra}. These  parameters are crucial in order to investigate the accuracy of interpolation  methods.

Considering a function $f$ to be recovered, it is possible, from a theoretical point of view, to take in account several combinations of function values and derivatives of $f$ and, moreover, this combinations may change from point to point of the set ${\cal X}$. In fact, the interpolant $H(u)$ in \eqref{bd} is such as to offer a complete flexibility.  On the other hand, in practice, the most interesting situation is when the values of $f$ and its first (and possibly second) derivatives are known at each point of ${\cal X}$. Hence, our numerical tests are restricted to this case and also to $\mathbb{S}^2$, as it is usually done using other interpolation methods, achieving  the advantage of permitting a comparison among different schemes (see e.g. \cite{fass99, fass07}). 

The test functions to be interpolated are taken from the restriction to $\mathbb{S}^2$ of the following trivariate functions
\begin{align*} 
f_1(x,y,z)=\frac{1}{10}[\exp x+ 2\exp(y+z)],\qquad f_2(x,y,z)=\sin x \sin y \sin z.
\end{align*}
Since the performance of the interpolant does not change significantly using other test functions (see e.g. \cite{hubbert01, pottmann90,  fass95}), for shortness we here report only the numerical results obtained considering $f_1$ and $f_2$. In Tables \ref{tab_1}--\ref{tab_2} we report the errors obtained for the interpolant \eqref{ch} using a complete Taylor expansion up to order zero (T0), one (T1) and two (T2). From these tables we can observe the significant improvement (in terms of accuracy) of the interpolant \eqref{ch} when making use of first and second derivatives in the Taylor expansion. Finally, to give an idea, in Table \ref{tab_lac} we show the results obtained in case of lacunary data, that is when a half of the first and second derivatives respectively are missing. This study points out that the interpolation scheme results in an unavoidable loss of accuracy due to the lack of information, but in any case the method can be applicable.  

\begin{table}[ht!]
{\small
		\begin{center}
			\begin{tabular}{|c|c|c|c|c|c|c|} \hline
			  & \multicolumn{2}{c|}{  \rule[-2mm]{0mm}{7mm}  T0} & \multicolumn{2}{c|}{  \rule[-2mm]{0mm}{7mm}  T1} & \multicolumn{2}{c|}{  \rule[-2mm]{0mm}{7mm}  T2} \\
			  \cline{2-7} \rule[-2mm]{0mm}{7mm}
			  $n$	 & MAE & RMSE & MAE & RMSE & MAE & RMSE \\
				\hline 
				\rule[0mm]{0mm}{3ex}
				$\hskip-2pt 500$   & $2.87{\rm E}-2$  &   $6.56{\rm E}-3$    & $2.96{\rm E}-3$  &   $1.12{\rm E}-3$	   & $1.89{\rm E}-4$  &   $2.38{\rm E}-5$    \\
			  \rule[0mm]{0mm}{3ex}
				$\hskip-2pt 1000$   & $2.07{\rm E}-2$  &   $4.09{\rm E}-3$    & $1.45{\rm E}-3$  &   $5.44{\rm E}-4$	   & $5.87{\rm E}-5$  &   $6.83{\rm E}-6$    \\
			  \rule[0mm]{0mm}{3ex}
				$\hskip-2pt 2000$   & $1.14{\rm E}-2$  &   $2.83{\rm E}-3$    & $5.88{\rm E}-4$  &   $2.67{\rm E}-4$	   & $1.73{\rm E}-5$  &   $2.24{\rm E}-6$    \\
			  \rule[0mm]{0mm}{3ex}
				$\hskip-2pt 4000$   & $8.41{\rm E}-3$  &   $1.94{\rm E}-3$    & $3.43{\rm E}-4$  &   $1.34{\rm E}-4$	   & $6.14{\rm E}-6$  &   $7.27{\rm E}-7$    \\
			  \rule[0mm]{0mm}{3ex}
				$\hskip-2pt 8000$   & $6.30{\rm E}-3$  &   $1.35{\rm E}-3$    & $1.83{\rm E}-4$  &   $6.65{\rm E}-5$	   & $1.98{\rm E}-6$  &   $2.25{\rm E}-7$    \\
			  \rule[0mm]{0mm}{3ex}
				$\hskip-2pt 16000$   & $4.83{\rm E}-3$  &   $9.08{\rm E}-4$    & $7.70{\rm E}-5$  &   $3.28{\rm E}-5$	   & $5.66{\rm E}-7$  &   $7.66{\rm E}-8$    \\
			  \hline 
			\end{tabular}
		\end{center}
		}
			\caption{MAEs and RMSEs computed on the sphere for $f_1$.}
			\label{tab_1}
	\end{table}

\begin{table}[ht!]
{\small
		\begin{center}
			\begin{tabular}{|c|c|c|c|c|c|c|} \hline
			  & \multicolumn{2}{c|}{  \rule[-2mm]{0mm}{7mm}  T0} & \multicolumn{2}{c|}{  \rule[-2mm]{0mm}{7mm}  T1} & \multicolumn{2}{c|}{  \rule[-2mm]{0mm}{7mm}  T2} \\
			  \cline{2-7} \rule[-2mm]{0mm}{7mm}
			  $n$	 & MAE & RMSE & MAE & RMSE & MAE & RMSE \\
				\hline 
				\rule[0mm]{0mm}{3ex}
				$\hskip-2pt 500$   & $1.44{\rm E}-2$  &   $3.72{\rm E}-3$    & $4.98{\rm E}-3$  &   $1.30{\rm E}-3$	   & $2.00{\rm E}-4$  &   $3.29{\rm E}-5$    \\
			  \rule[0mm]{0mm}{3ex}
				$\hskip-2pt 1000$   & $8.47{\rm E}-3$  &   $2.37{\rm E}-3$    & $2.90{\rm E}-3$  &   $6.36{\rm E}-4$	   & $4.09{\rm E}-5$  &   $1.08{\rm E}-5$    \\
			  \rule[0mm]{0mm}{3ex}
				$\hskip-2pt 2000$   & $6.78{\rm E}-3$  &   $1.70{\rm E}-3$    & $1.22{\rm E}-3$  &   $3.13{\rm E}-4$	   & $1.72{\rm E}-5$  &   $3.76{\rm E}-6$    \\
			  \rule[0mm]{0mm}{3ex}
				$\hskip-2pt 4000$   & $5.16{\rm E}-3$  &   $1.20{\rm E}-3$    & $6.70{\rm E}-4$  &   $1.55{\rm E}-4$	   & $6.35{\rm E}-6$  &   $1.28{\rm E}-6$    \\
			  \rule[0mm]{0mm}{3ex}
				$\hskip-2pt 8000$   & $3.36{\rm E}-3$  &   $8.70{\rm E}-4$    & $2.52{\rm E}-4$  &   $7.66{\rm E}-5$	   & $2.16{\rm E}-6$  &   $4.60{\rm E}-7$    \\
			  \rule[0mm]{0mm}{3ex}
				$\hskip-2pt 16000$   & $3.35{\rm E}-3$  &   $6.20{\rm E}-4$    & $1.15{\rm E}-4$  &   $3.78{\rm E}-5$	   & $5.48{\rm E}-7$  &   $1.57{\rm E}-7$    \\
			  \hline 
			\end{tabular}
		\end{center}
		}
			\caption{MAEs and RMSEs computed on the sphere for $f_2$.}
			\label{tab_2}
	\end{table}

\begin{table}[ht!]
{\small
		\begin{center}
			\begin{tabular}{|c|c|c|c|c|} \hline
			  & \multicolumn{2}{c|}{  \rule[-2mm]{0mm}{7mm}  missing 1st der.} & \multicolumn{2}{c|}{  \rule[-2mm]{0mm}{7mm}  missing 2nd der.}  \\
			  \cline{2-5} \rule[-2mm]{0mm}{7mm}
			  $n$	 & MAE & RMSE & MAE & RMSE  \\
				\hline 
				\rule[0mm]{0mm}{3ex}
				$\hskip-2pt 500$   & $3.16{\rm E}-2$  &   $5.68{\rm E}-3$    & $3.25{\rm E}-3$  &   $9.17{\rm E}-4$	 \\
			  \rule[0mm]{0mm}{3ex}
				$\hskip-2pt 1000$   & $2.20{\rm E}-2$  &   $3.35{\rm E}-3$    & $1.80{\rm E}-3$  &   $4.47{\rm E}-4$	  \\
			  \rule[0mm]{0mm}{3ex}
				$\hskip-2pt 2000$   & $1.19{\rm E}-2$  &   $2.13{\rm E}-3$    & $7.86{\rm E}-4$  &   $2.23{\rm E}-4$	  \\
			  \rule[0mm]{0mm}{3ex}
				$\hskip-2pt 4000$   & $8.28{\rm E}-3$  &   $1.43{\rm E}-3$    & $4.02{\rm E}-4$  &   $1.13{\rm E}-4$	  \\
			  \rule[0mm]{0mm}{3ex}
				$\hskip-2pt 8000$   & $6.40{\rm E}-3$  &   $9.38{\rm E}-4$    & $1.90{\rm E}-4$  &   $5.73{\rm E}-5$	  \\
			  \rule[0mm]{0mm}{3ex}
				$\hskip-2pt 16000$   & $4.90{\rm E}-3$  &   $6.63{\rm E}-4$    & $9.30{\rm E}-5$  &   $2.84{\rm E}-5$	  \\
			  \hline 
			\end{tabular}
		\end{center}
			\caption{MAEs and RMSEs computed on the sphere with lacunary data for $f_1$. Left: missing a half of the 1st derivatives, right: missing a half of the 2nd derivatives.}
			\label{tab_lac}
			}
	\end{table}

The considered interpolation scheme for the sphere is  suitable for parallel implementation (see \cite{allasia97, costanzo05, derossi07b}). In the implementation of the parallel algorithm for a distributed memory machine, the data are assigned to  $p$ processors by breaking the set ${\cal X}$ into subsets ${\cal X}_k,~k=1, \ldots, p$. In this way each processor proceeds to solve the interpolation problem on a subset ${\cal X}_k$. The parallel algorithm consists of three steps: a) partitioning and data distribution, so that each subset has an approximately equal number of points, b) local interpolation solving, after having determined the radius of influence for each point of the considered subset, c) data collection and evaluation phase, where each slave processor sends its partial results to the master processor. In the ideal case, when the algorithm is completely and efficiently parallelizable, the speed-up must be equal to the number of processors involved.


\section{Numerical Results on Riemannian Manifolds}

Moving from the considerations on the sphere to those on  Riemannian manifolds in general,  nothing  changes for what concerns the structure of the interpolant \eqref{ch} and of the cardinal basis functions to be used. Instead, what changes dramatically is the problem of calculating the geodesic distance between points, because on the sphere one has got a simple analytic expression of the geodesic distance  while this does not happen for other manifolds.

\subsection{Computation of the Geodesic Distance on Manifolds}

  In order to study properties of geodesics on a $m$-dimensional Riemannian manifold ${\cal M}$, it is convenient to consider a  connected chart $(U, \varphi)$ and for each $u\in U$ the vector $\varphi (u)=(v_1, \ldots, v_m)\in V$ of local coordinates. Then,  a geodesic   on $U$ is a curve given  by $m$ functions $v_1(s), \ldots, v_m(s)$ which satisfy the system of second order differential equations, called {\sl geodesic equations},
\begin{equation}
\frac{d^2v_k}{ds^2}+\sum_{i,j=1}^m\Gamma_{ij}^k
\frac{dv_i}{ds}
\frac{dv_j}{ds}
=0, \quad k=1, \ldots, m,
\label{fa}
\end{equation} 
where   $s$ is the arclength parameter and $\Gamma_{ij}^k$ are the Christoffel symbols of the second kind. It is possible to express the Christoffel simbols in terms of the components  $g_{ks}$ of the metric matrix  and their derivatives as follows
\begin{equation}
\Gamma _{ij}^k=\frac{1}{2}\sum_{s=1}^m g^{ks}\Big(\frac{\partial g_{si}}{\partial v_j}+\frac{\partial g_{js}}{\partial v_i}-\frac{\partial g_{ij}}{\partial v_s}\Big),
\label{faa}
\end{equation}
where  $g^{ks}$ are the  components  of the matrix $(g^{ij})$, inverse of the matrix $(g_{ij})$.
As the manifold has dimension $m$, the geodesic equations are a system of $m$ ordinary differential equations for the $m$  variables $v_k$. Thus, allied with initial conditions consisting of a point on the manifold and a tangent vector at the point, the system can  theoretically be uniquely solved, at least locally, but actually  there are  serious difficulties.

Computing the geodesic distance is less prohibitive if we merely consider a particular, but important, subclass of Riemannian manifolds, namely the regular surfaces in $\RR^3$ parametrically defined on an open set of $\RR^2$ 
 by a map of the type $u=(u_1(v_1,  v_2), u_2(v_1, v_2),$  $u_3(v_1,  v_2))$.
In this case  the  system \eqref{fa} reduces to two equations.

A further simplification is achieved considering  for the  set $U$ on the surface ${\cal M}$ orthogonal local coordinates, so that  $g_{12}=g_{21}=0$. For  $\gamma=\gamma(s)$ to be a geodesic on  $U$, then it is necessary and sufficient that the {geodesic equations}
\begin{align}
\frac{d^2 v_1}{d s^2}+\frac{1}{2g_{11}}\frac{\partial g_{11}}{\partial v_1}\Big(\frac{d v_1}{d s}\Big)^2+\frac{1}{g_{11}}\frac{\partial g_{11}}{\partial v_2}\frac{d v_1}{d s}\frac{d v_2}{d s}-\frac{1}{2g_{11}}\frac{\partial g_{22}}{\partial v_1}\Big(\frac{d v_2}{d s}\Big)^2&=0,\cr
\frac{d^2 v_2}{d s^2}-\frac{1}{2g_{22}}\frac{\partial g_{11}}{\partial v_2}\Big(\frac{d v_1}{d s}\Big)^2+\frac{1}{g_{22}}\frac{\partial g_{22}}{\partial v_1}\frac{d v_1}{d s}\frac{d v_2}{d s}+\frac{1}{2g_{22}}\frac{\partial g_{22}}{\partial v_2}\Big(\frac{d v_2}{d s}\Big)^2&=0,
\label{fb}
\end{align}
are satisfied, expressing now the Christoffel symbols  by \eqref{faa}.

An even more favorable situation  is achieved  considering the Clairaut parametrizations. We say that an orthogonal  parametrization is a {\sl Clairaut parametrization in} $v_1$ if 
$$
\frac{\partial g_{11}}{\partial v_2}=\frac{\partial g_{22}}{\partial v_2}=0.
$$
Similarly,  we say that an orthogonal  parametrization is a {\sl Clairaut parametrization in} $v_2$ if 
$$
\frac{\partial g_{11}}{\partial v_1}=\frac{\partial g_{22}}{\partial v_1}=0.
$$
The geodesic equations simplify in these cases to the $v_1$-{\sl Clairaut geodesic equations}
\begin{align}
\frac{d^2 v_1}{d s^2}+\frac{1}{2g_{11}}\frac{\partial g_{11}}{\partial v_1}\Big(\frac{d v_1}{d s}\Big)^2
-\frac{1}{2g_{11}}\frac{\partial g_{22}}{\partial v_1}\Big(\frac{d v_2}{d s}\Big)^2&=0,\cr
\frac{d^2 v_2}{d s^2}
+\frac{1}{g_{22}}\frac{\partial g_{22}}{\partial v_1}\frac{d v_1}{d s}\frac{d v_2}{d s}
&=0,
\label{fc}
\end{align}
and
to the $v_2$-{\sl Clairaut geodesic equations}
\begin{align}
\frac{d^2 v_1}{d s^2}
+\frac{1}{g_{11}}\frac{\partial g_{11}}{\partial v_2}\frac{d v_1}{d s}\frac{d v_2}{d s}
&
=0,\cr
\frac{d^2 v_2}{d s^2}
-\frac{1}{2g_{22}}\frac{\partial g_{11}}{\partial v_2}\Big(\frac{d v_1}{d s}\Big)^2
+\frac{1}{2g_{22}}\frac{\partial g_{22}}{\partial v_2}\Big(\frac{d v_2}{d s}\Big)^2
&=0,
\label{fd}
\end{align}
respectively.
As an example, 
considering the torus
\begin{equation}
{\cal T}^2(v_1, v_2)=\{((R+r\cos v_1)\cos v_2,~(R+r\cos v_1)\sin v_2, ~r\sin v_1 ): v_1, v_2\in [0, 2\pi[ \}
\label{fh}
\end{equation}
the  $v_1$-{Clairaut geodesic equations}  are
\begin{align}
\frac{d^2 v_1}{d s^2}+
\frac{R+r\cos v_1}{r}\sin v_1\Big(\frac{dv_2}{ds}\Big)^2&=0,\cr
\frac{d^2 v_2}{d s^2}-2\frac{r\sin v_1}{R+r\cos v_1}\frac{dv_1}{ds}\frac{dv_2}{ds}
&=0.
\label{fk}
\end{align}

Actually  the surfaces of revolution, which include many  cases important for applications,  appear to be the most  manageable. 
Without loss of generalization, let us consider a plane $\pi\subset\RR^3$ generated by the unit vectors $e_1, e_3\in\RR^3$, a straight line $l$ generated by the unit vector $e_3$, and a curve $c\in C^k, ~k\ge 1$, which is disjoint from $l$ and belongs to the positive halfplane with respect to $e_1$.  Rotating the curve $c$ around the line $l$ we obtain a  surface of revolution ${\cal M}$ with generating curve  $c$ and revolution axis $l$. If $c$ is parametrically represented  by
$$
c(v_1)=
\{(\alpha(v_1),  \beta(v_1)): v_1\in[a, b]\subset\RR,~ \alpha(v_1)>0 \},
$$
the surface is given by
$$
{\cal M}(v_1, v_2)=\{(\alpha(v_1)\cos v_2, \alpha(v_1)\sin v_2, \beta(v_1))\in\RR^3: v_1\in[a, b],~ v_2\in [0, 2\pi[,~\alpha(v_1)>0\}.
$$
For a fixed $v_1=v_1^0$ the curve $p={\cal M}(v_1^0, v_2)\subset\RR^3$
is called a {\sl parallel} of ${\cal M}$ and represents the circle of radius $\alpha (v_1^0) $ obtained by rotating the point $(\alpha(v_1^0), \beta(v_1^0))\in c$ around the line $l$. 
Similarly, for a fixed $v_2=v_2^0$ the curve $m={\cal M}(v_1, v_2^0)\subset\RR^3$
is called a {\sl meridian} of ${\cal M}$ and  is obtained by rotating $c$ of an angle  $v_2^0$ around  $l$.

If the curve is parameterized with respect to the arclength $s$, the differential equations of geodesics for surfaces of rotation are
\begin{align*}
&\frac{d^2 v_1}{ds^2}-\alpha\frac{d \alpha}{ds}\Big(\frac{dv_2}{ds}\Big)^2=0, \cr
&\frac{d^2v_2}
{ds^2}+
\frac{2}{\alpha}
\frac{d\alpha}{ds}
\frac{dv_1}{ds}\frac{dv_2}{ds}=0. 
 \end{align*} 
Important consequences of these equations are:
\begin{enumerate}
\item[i)] the meridians of a surface of revolution are geodesic curves,
\item[ii)] a parallel is a geodesic curve if and only if it is obtained by rotating a point on the generating curve whose tangent vector is parallel to the axis of revolution.
\end{enumerate}
Note that the geodesic equations of surfaces of revolution parameterised with respect to the arclength are  particular cases of $v_1$-parameterization of Clairaut.
 
We used \textsc{Matlab} to numerically solve equations of geodesics on a parametric surface and to picture the relevant graphs. If ${\cal M}\subset\RR^3$ is the surface 
$$
\varphi^{-1}(v_1, v_2)=\{(u_1, u_2, u_3):(u_1, u_2, u_3)\in U\subset {\cal M}, (v_1, v_1)\in \varphi (U)\subset\RR^2\},
$$
we built two programs, which are similar but meet different needs. One of them  resolves the system \eqref{fa} with $m=2$ (as well as its particular cases) and finds $(v_1(s), v_2(s))$ with the arclength parameter $s$ or a multiple of it, starting from the initial point $(v_1 (0), v_2 (0))$ and the derivatives $(dv_1(0)/ ds, dv_2(0)/ ds)$.  Then the program draws the support of ${\cal M}$ by varying $(v_1, v_2)\in \varphi(U)$  and traces the geodesic curve leaving  $s$ to vary in a given  interval $[s_i, s_e]\subset\RR$ required in input. Seeing pictures of geodesics is obviously interesting and useful, but it is not our primary goal (see e.g. \cite{oprea07,abdelall13}). 

The other program resolves the system \eqref{fa} with $m=2$ (as well as its particular cases) and finds $(v_1(s), v_2(s))$ with the arclength parameter $s$ or a multiple of it, given the initial point $(v_1 (s_i), v_2 (s_i))$ and the end point  $(v_1(s_e), v_2(s_e))$. Then, the program traces the geodesic connecting the two points and, above all, compute the geodesic distance between them. This second program starts with an approximate path of the geodesic and improves the solution iteratively. Since the geodesic distance between the initial and the end points is very small, we can choose a segment as the initial guess. In general, the method works well, but in a few cases the convergence is not assured despite requiring compactness (see e.g. \cite{maekawa96,kimmel98}).

\subsection{Computation of the Interpolant on Manifolds}

To optimize the nearest neighbour searching procedure for data on a general Riemannian manifold, one can continue to use the  techniques already described for the sphere. Of course, in individual cases, more difficulties than with the sphere may arise in the use of those procedures. 

To test the performance of our interpolation method, we need to get a uniform (or quasi-uniform) distribution of nodes on the considered Riemannian manifolds. Unfortunately, finding a convenient distribution is another critical point. To face the problem of generating a uniform distribution of points on analytic surfaces, there appear to be interesting the results of some recent papers \cite{kopytov15,kopytov12,kopytov15a,melfi04,petrillo05}. On the other hand, to  get information on uniformity and density  of the  distribution of points in the set ${\cal X}\subset U$, the {separation distance} \eqref{cdb} and the {fill distance} \eqref{cra}, already considered for the sphere, continue to be crucial parameters to assess the accuracy of interpolation methods.

\begin{figure}[ht!]
		\centering
		\includegraphics[height=.31\textheight]{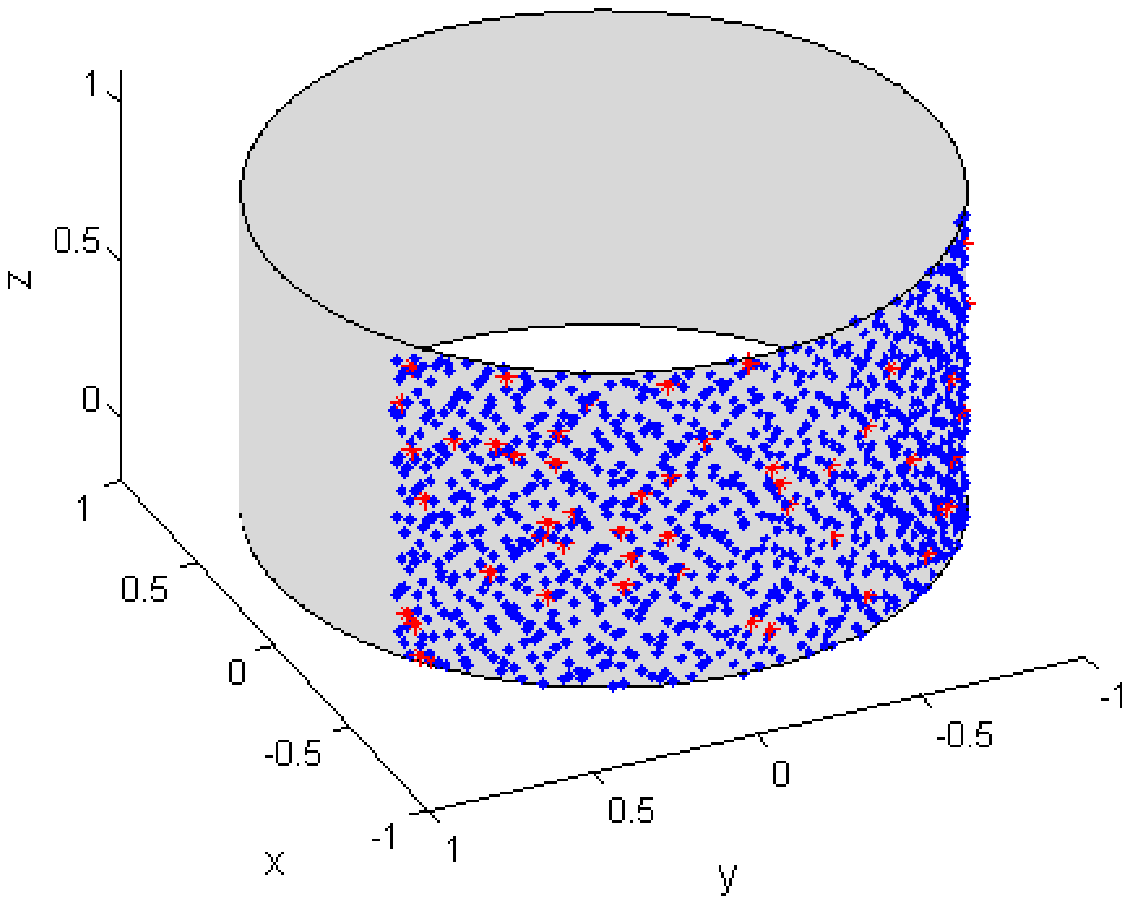} \hskip -1.8cm
		\includegraphics[height=.31\textheight]{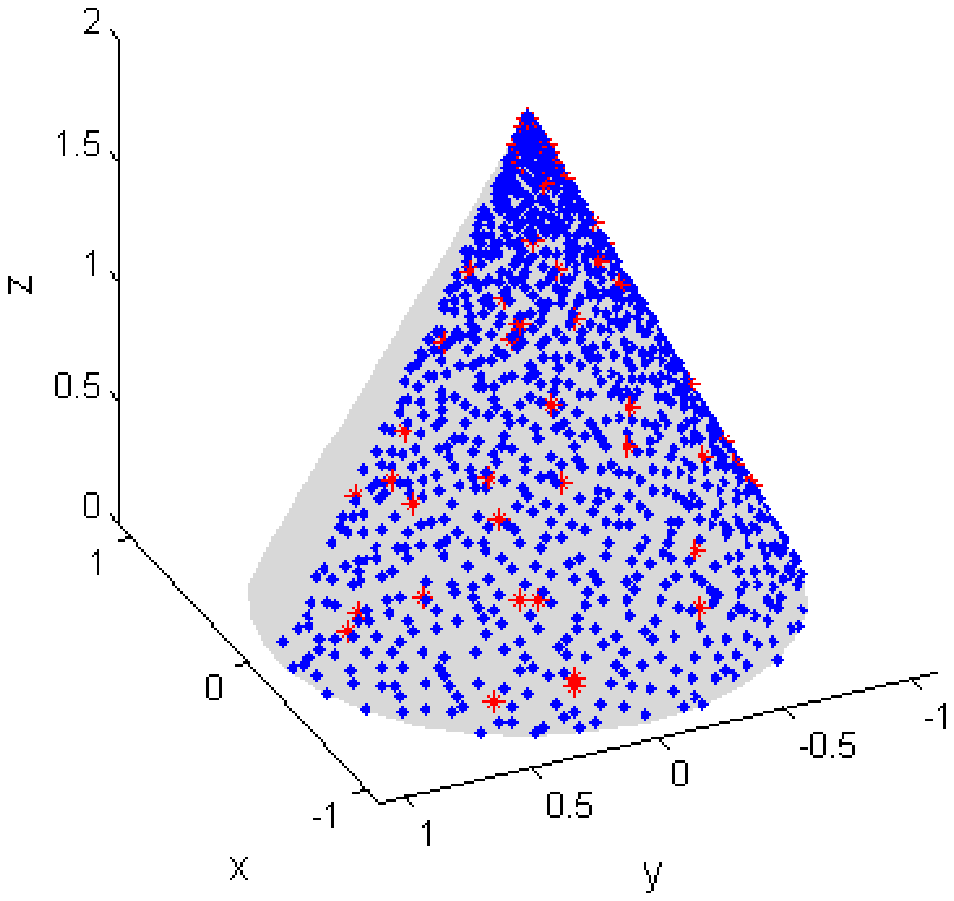}
		\caption{Example of interpolation (blue dot) and evaluation (red star) points on the cylinder (left) and the cone (right).}
		\label{fig_cylinder_cone}
	\end{figure}

To test our interpolant on manifolds, we focus on cylinder and cone. As interpolation nodes, we take some sets of $n$ uniformly random Halton data points, originally contained in the unit square $[0, 1]^2 \subset \mathbb{R}^2$ and then mapped onto the surface of cylinder and cone via, respectively, the equations
\begin{align*}
x = r\cos(2\pi p), \quad y = r\sin(2\pi p), \quad z = q, 
\end{align*} 
and
\begin{align*}
x = \frac{(h-z)}{h}r\cos(2\pi p), \quad y = \frac{(h-z)}{h}r\sin(2\pi p), \quad z = hq, 
\end{align*}
where $(p,q) \in [0,1]^2$, $r$ is the radius and $h$ denotes the height of the cone. The computation of interpolation errors, using \eqref{MAE_RMSE}, is carried out mapping as earlier a set of $n_{eval} = 50$ evaluation points generated by the \texttt{rand} command of \textsc{Matlab}. Note that in our tests we consider a chart of ${\cal M}$, taking all points belonging to the cylinder for $x<-0.5$ and to the cone for $x<0$, assuming $r=1$ and $h=2$. An example of interpolation and evaluation points defined on the charts of cylinder and cone is given in Figure \ref{fig_cylinder_cone}.

In order to recover a function $f$,  known on a set ${\cal X}\subset U\subset{\cal M}$ together with some of its derivatives,  the interpolant $H(u)$ in \eqref{bd} is quite efficient for any combination of derivatives.  However, actually,  the most interesting situation is when the values of $f$ and its first and second derivatives are known at each point of ${\cal X}$. Hence, our numerical tests are restricted to this case and also to a   chart  of the cylinder and a chart of the cone. As regard to Hermite-Birkhoff interpolation on manifolds  it is difficult to find  numerical tests in the literature, as far as we know, while theoretical considerations are not lacking (see e.g. \cite{dyn97, dyn99, narc95}).
The test functions to be interpolated are taken from the restriction to $U\subset {\cal M}$ of the trivariate functions already considered for the sphere, but we report only the numerical results obtained considering $f_1$ and $f_2$. Interpolation errors computed for the interpolant \eqref{ch} using a complete Taylor expansion up to order zero (T0), one (T1) and two (T2) are shown in Tables \ref{tab_3}--\ref{tab_4} for the cylinder and Tables \ref{tab_5}--\ref{tab_6} for the cone. From these numerical experiments we obtain an error behavior similar to that observed in the previous section for the sphere.

\begin{table}[ht!]
{\small
		\begin{center}
			\begin{tabular}{|c|c|c|c|c|c|c|} \hline
			  & \multicolumn{2}{c|}{  \rule[-2mm]{0mm}{7mm}  T0} & \multicolumn{2}{c|}{  \rule[-2mm]{0mm}{7mm}  T1} & \multicolumn{2}{c|}{  \rule[-2mm]{0mm}{7mm}  T2} \\
			  \cline{2-7} \rule[-2mm]{0mm}{7mm}
			  $n$	 & MAE & RMSE & MAE & RMSE & MAE & RMSE \\
				\hline 
				\rule[0mm]{0mm}{3ex}
$\hskip-2pt 500$   & $4.01{\rm E}-2$  &   $1.15{\rm E}-2$    & $4.22{\rm E}-3$  &   $1.36{\rm E}-3$	   & $1.38{\rm E}-4$  &   $3.18{\rm E}-5$    \\
			  \rule[0mm]{0mm}{3ex}
$\hskip-2pt 1000$   & $2.37{\rm E}-2$  &   $6.19{\rm E}-3$    & $1.29{\rm E}-3$  &   $5.22{\rm E}-4$	   & $1.69{\rm E}-5$  &   $5.67{\rm E}-6$    \\
			  \rule[0mm]{0mm}{3ex}
$\hskip-2pt 2000$   & $1.27{\rm E}-2$  &   $3.63{\rm E}-3$    & $8.71{\rm E}-4$  &   $2.99{\rm E}-4$	   & $1.06{\rm E}-5$  &   $2.48{\rm E}-6$    \\
			  \rule[0mm]{0mm}{3ex}
$\hskip-2pt 4000$   & $5.09{\rm E}-3$  &   $1.87{\rm E}-3$    & $3.66{\rm E}-4$  &   $1.34{\rm E}-4$	   & $2.74{\rm E}-6$  &   $6.79{\rm E}-7$    \\
			  \rule[0mm]{0mm}{3ex}
$\hskip-2pt 8000$   & $4.74{\rm E}-3$  &   $1.19{\rm E}-3$    & $1.27{\rm E}-4$  &   $6.56{\rm E}-5$	   & $4.79{\rm E}-7$  &   $1.44{\rm E}-7$    \\
			  \rule[0mm]{0mm}{3ex}
$\hskip-2pt 16000$  & $2.93{\rm E}-3$  &   $1.07{\rm E}-3$    & $7.87{\rm E}-5$  &   $3.55{\rm E}-5$	   & $1.82{\rm E}-7$  &   $6.81{\rm E}-8$  \\
			  \hline 
			\end{tabular}
		\end{center}
		}
			\caption{MAEs and RMSEs computed on the cylinder for $f_1$.}
			\label{tab_3}
	\end{table}

\begin{table}[ht!]
{\small
		\begin{center}
			\begin{tabular}{|c|c|c|c|c|c|c|} \hline
			  & \multicolumn{2}{c|}{  \rule[-2mm]{0mm}{7mm}  T0} & \multicolumn{2}{c|}{  \rule[-2mm]{0mm}{7mm}  T1} & \multicolumn{2}{c|}{  \rule[-2mm]{0mm}{7mm}  T2} \\
			  \cline{2-7} \rule[-2mm]{0mm}{7mm}
			  $n$	 & MAE & RMSE & MAE & RMSE & MAE & RMSE \\
				\hline 
				\rule[0mm]{0mm}{3ex}
$\hskip-2pt 500$   & $1.91{\rm E}-2$  &   $5.13{\rm E}-3$    & $2.59{\rm E}-3$  &   $9.23{\rm E}-4$	   & $7.94{\rm E}-5$  &   $2.44{\rm E}-5$    \\
			  \rule[0mm]{0mm}{3ex}
$\hskip-2pt 1000$   & $6.16{\rm E}-3$  &   $2.85{\rm E}-3$    & $1.09{\rm E}-3$  &   $3.85{\rm E}-4$	   & $1.74{\rm E}-5$  &   $6.65{\rm E}-6$    \\
			  \rule[0mm]{0mm}{3ex}
$\hskip-2pt 2000$   & $3.96{\rm E}-3$  &   $1.55{\rm E}-3$    & $4.99{\rm E}-4$  &   $1.79{\rm E}-4$	   & $5.41{\rm E}-6$  &   $1.75{\rm E}-6$    \\
			  \rule[0mm]{0mm}{3ex}
$\hskip-2pt 4000$   & $2.73{\rm E}-3$  &   $1.10{\rm E}-3$    & $2.21{\rm E}-4$  &   $9.47{\rm E}-5$	   & $1.59{\rm E}-6$  &   $6.09{\rm E}-7$    \\
			  \rule[0mm]{0mm}{3ex}
$\hskip-2pt 8000$   & $2.15{\rm E}-3$  &   $7.33{\rm E}-4$    & $8.83{\rm E}-5$  &   $4.28{\rm E}-5$	   & $5.43{\rm E}-7$  &   $1.83{\rm E}-7$    \\
			  \rule[0mm]{0mm}{3ex}
$\hskip-2pt 16000$  & $1.23{\rm E}-3$  &   $5.20{\rm E}-4$    & $5.66{\rm E}-5$  &   $2.51{\rm E}-5$	   & $1.93{\rm E}-7$  &   $6.33{\rm E}-8$  \\
			  \hline 
			\end{tabular}
		\end{center}
		}
			\caption{MAEs and RMSEs computed on the cylinder for $f_2$.}
			\label{tab_4}
	\end{table}

\begin{table}[ht!]
{\small
		\begin{center}
			\begin{tabular}{|c|c|c|c|c|c|c|} \hline
			  & \multicolumn{2}{c|}{  \rule[-2mm]{0mm}{7mm}  T0} & \multicolumn{2}{c|}{  \rule[-2mm]{0mm}{7mm}  T1} & \multicolumn{2}{c|}{  \rule[-2mm]{0mm}{7mm}  T2} \\
			  \cline{2-7} \rule[-2mm]{0mm}{7mm}
			  $n$	 & MAE & RMSE & MAE & RMSE & MAE & RMSE \\
				\hline 
				\rule[0mm]{0mm}{3ex}
$\hskip-2pt 500$   & $2.29{\rm E}-2$  &   $9.64{\rm E}-3$    & $3.20{\rm E}-3$  &   $1.62{\rm E}-3$	   & $7.80{\rm E}-5$  &   $2.41{\rm E}-5$    \\
			  \rule[0mm]{0mm}{3ex}
$\hskip-2pt 1000$   & $1.16{\rm E}-2$  &   $5.64{\rm E}-3$    & $1.52{\rm E}-3$  &   $7.35{\rm E}-4$	   & $3.95{\rm E}-5$  &   $9.93{\rm E}-6$    \\
			  \rule[0mm]{0mm}{3ex}
$\hskip-2pt 2000$   & $6.60{\rm E}-3$  &   $3.17{\rm E}-3$    & $1.01{\rm E}-3$  &   $3.76{\rm E}-4$	   & $1.06{\rm E}-5$  &   $2.97{\rm E}-6$    \\
			  \rule[0mm]{0mm}{3ex}
$\hskip-2pt 4000$   & $7.81{\rm E}-3$  &   $2.51{\rm E}-3$    & $3.50{\rm E}-4$  &   $1.72{\rm E}-4$	   & $3.89{\rm E}-6$  &   $9.76{\rm E}-7$    \\
			  \rule[0mm]{0mm}{3ex}
$\hskip-2pt 8000$   & $4.11{\rm E}-3$  &   $1.47{\rm E}-3$    & $1.68{\rm E}-4$  &   $8.78{\rm E}-5$	   & $8.80{\rm E}-7$  &   $2.64{\rm E}-7$    \\
			  \rule[0mm]{0mm}{3ex}
$\hskip-2pt 16000$  & $2.74{\rm E}-3$  &   $9.02{\rm E}-4$    & $8.58{\rm E}-5$  &   $4.50{\rm E}-5$	   & $2.43{\rm E}-7$  &   $8.61{\rm E}-8$  \\
			  \hline 
			\end{tabular}
		\end{center}
		}
			\caption{MAEs and RMSEs computed on the cone for $f_1$.}
			\label{tab_5}
	\end{table}

\begin{table}[ht!]
{\small
		\begin{center}
			\begin{tabular}{|c|c|c|c|c|c|c|} \hline
			  & \multicolumn{2}{c|}{  \rule[-2mm]{0mm}{7mm}  T0} & \multicolumn{2}{c|}{  \rule[-2mm]{0mm}{7mm}  T1} & \multicolumn{2}{c|}{  \rule[-2mm]{0mm}{7mm}  T2} \\
			  \cline{2-7} \rule[-2mm]{0mm}{7mm}
			  $n$	 & MAE & RMSE & MAE & RMSE & MAE & RMSE \\
				\hline 
				\rule[0mm]{0mm}{3ex}
$\hskip-2pt 500$   & $1.18{\rm E}-2$  &   $3.90{\rm E}-3$    & $3.52{\rm E}-3$  &   $1.35{\rm E}-3$	   & $1.58{\rm E}-4$  &   $4.32{\rm E}-5$    \\
			  \rule[0mm]{0mm}{3ex}
$\hskip-2pt 1000$   & $9.06{\rm E}-3$  &   $2.39{\rm E}-3$    & $1.93{\rm E}-3$  &   $7.41{\rm E}-4$	   & $4.31{\rm E}-5$  &   $1.29{\rm E}-5$    \\
			  \rule[0mm]{0mm}{3ex}
$\hskip-2pt 2000$   & $4.72{\rm E}-3$  &   $1.31{\rm E}-3$    & $8.12{\rm E}-4$  &   $3.24{\rm E}-4$	   & $1.57{\rm E}-5$  &   $5.55{\rm E}-6$    \\
			  \rule[0mm]{0mm}{3ex}
$\hskip-2pt 4000$   & $3.47{\rm E}-3$  &   $9.21{\rm E}-4$    & $3.36{\rm E}-4$  &   $1.61{\rm E}-4$	   & $6.62{\rm E}-6$  &   $1.59{\rm E}-6$    \\
			  \rule[0mm]{0mm}{3ex}
$\hskip-2pt 8000$   & $2.29{\rm E}-3$  &   $6.63{\rm E}-4$    & $2.11{\rm E}-4$  &   $9.30{\rm E}-5$	   & $1.85{\rm E}-6$  &   $4.63{\rm E}-7$    \\
			  \rule[0mm]{0mm}{3ex}
$\hskip-2pt 16000$  & $1.10{\rm E}-3$  &   $3.92{\rm E}-4$    & $1.22{\rm E}-4$  &   $4.40{\rm E}-5$	   & $5.27{\rm E}-7$  &   $1.56{\rm E}-7$  \\
			  \hline 
			\end{tabular}
		\end{center}
		}
			\caption{MAEs and RMSEs computed on the cone for $f_2$.}
			\label{tab_6}
	\end{table}

The considered interpolation schemes for the Riemannian manifolds are suitable for parallel implementation as explained for the sphere.


	
\section*{Acknowledgements}		
		
This work was supported by the University of Turin via grant \lq\lq Metodi numerici nelle scienze applicate\rq\rq.

\end{document}